\documentclass[11pt]{article}

\usepackage[T1]{fontenc}
\usepackage{lmodern}
\usepackage{amsmath,amssymb,amsthm}
\usepackage{graphicx}
\usepackage[margin=1.12in]{geometry}
\usepackage{microtype}
\usepackage{url}
\usepackage[hidelinks]{hyperref}
\newtheorem{theorem}{Theorem}
\newtheorem{lemma}[theorem]{Lemma}
\newtheorem{proposition}[theorem]{Proposition}
\newtheorem{corollary}[theorem]{Corollary}

\newtheorem*{theorem*}{Theorem}
\newtheorem*{corollary*}{Corollary}

\theoremstyle{remark}
\newtheorem{remark}[theorem]{Remark}

\theoremstyle{plain}

\title{Component Trees of Pachner Graphs\\
with Applications to Triangulated \(3\)-Spheres}
\author{Vance Faber\thanks{Hoquiam, WA, \texttt{vance.faber@gmail.com}}
\and Michael Murphy\thanks{Venice, CA}}
\date{}

\begin{document}
\maketitle

\begin{abstract}
Fix a closed connected PL \(d\)-manifold \(M\), where \(d\ge2\).  For each
\(n\), let \(G(n)\) be the graph whose vertices are the simplicial
triangulations of \(M\) with \(n\) vertices, with an edge for each
vertex-preserving bistellar flip, and let \(C(n)\) be its set of connected
components.  We form a directed graph \(\mathcal C\) whose vertices are the
elements of \(C(n)\), over all \(n\).  There is an edge from
\(U\in C(n)\) to \(V\in C(n+1)\) when a facet subdivision---equivalently, a
\(1\)--\((d+1)\) bistellar move---takes some triangulation in \(U\) to a
triangulation in \(V\).

We prove that every vertex of \(\mathcal C\) has out-degree one: every facet
subdivision of every triangulation in \(U\) lies in the same component of
\(G(n+1)\).  Since the full Pachner graph of \(M\) is connected, the
underlying undirected graph of \(\mathcal C\) is a graded tree.  Consequently,
the directed rays from any two components eventually meet.  The unique path
between two components, which we call their \emph{canonical component path},
lifts to a Pachner sequence consisting first of facet subdivisions, then of
vertex-preserving bistellar moves, and finally of inverse facet subdivisions.

A chosen triangulation determines a monotone closure containing exactly one
fixed-vertex component at each higher level.  Let \(n_0\) be the minimum
number of vertices in a triangulation of \(M\).  When \(G(n_0)\) is connected,
its stabilization ray is a distinguished \emph{trunk}.  For \(S^3\), the
trunk begins at the boundary of the \(4\)-simplex.  We prove that it contains
every triangulated \(3\)-sphere with at most \(11\) vertices, every
\(12\)-vertex seed sphere with minimum edge valence at least \(4\), and every
facet subdivision of each of the four known unflippable \(3\)-spheres.  For
four surfaces whose flip graphs are disconnected at the minimum level, all
components merge after one facet subdivision; one example merges \(59\)
isolated components at once.  Examples from the \(3\)-torus and the
Bagchi--Datta family illustrate unflippable behavior beyond \(S^3\) and in
arbitrarily high dimensions.
\end{abstract}

\section{Introduction}
\label{sec:introduction}

Pachner's theorem states that any two PL-homeomorphic triangulations of a
closed manifold are related by a finite sequence of bistellar moves
\cite{P91}.  Thus, the full Pachner graph of a fixed closed PL manifold is
connected.  Its horizontal slices, obtained by fixing the number of vertices
and allowing only vertex-preserving bistellar moves, need not be connected.
The purpose of this paper is to describe how the components of these slices
fit together as the number of vertices increases.

Fix a closed connected PL \(d\)-manifold \(M\), where \(d\ge2\).  Let
\(G(n)\) be the graph whose vertices are the simplicial triangulations of
\(M\) with \(n\) vertices, with an edge for each vertex-preserving bistellar
move.  Write
\[
 C(n)=\pi_0\bigl(G(n)\bigr)
\]
for the set of connected components of \(G(n)\).  Since \(M\) is fixed, we
suppress it from the notation.  When several manifolds are compared, we
restore the subscript and write \(G_M(n)\) and \(C_M(n)\).

Form a directed graph \(\mathcal C\) whose vertices are the elements of
\(C(n)\), over all \(n\).  There is a directed edge from \(U\in C(n)\) to
\(V\in C(n+1)\) if some facet subdivision of some triangulation in \(U\)
belongs to \(V\).  Our main theorem is the following.

\begin{theorem*}[Component Stabilization Theorem]
Every vertex of \(\mathcal C\) has out-degree one.  Equivalently, if \(P\)
and \(Q\) lie in the same component of \(G(n)\), then every facet subdivision
of \(P\) is connected by vertex-preserving bistellar moves to every facet
subdivision of \(Q\).
\end{theorem*}

The same local argument works in every dimension.  After a new vertex is
inserted into a facet, it can be slid across a common \((d-1)\)-face,
called a ridge, into an adjacent facet by a \(2\)--\(d\) move followed by a
\(d\)--\(2\) move.  Since the
facet-adjacency graph of a connected closed triangulated manifold is
connected, all choices of the subdivided facet give triangulations in one
component of \(G(n+1)\).  A vertex-preserving bistellar move downstairs can
then be performed upstairs by choosing the subdivided facet outside its
support.

The global consequences follow from the Component Stabilization Theorem
together with Pachner's theorem.  The full Pachner graph is connected, and
the grading of \(\mathcal C\), together with the unique upward edge from each
vertex, rules out cycles.  Hence the underlying undirected graph of
\(\mathcal C\) is a graded tree.  The directed rays from any two components
eventually meet, and the unique path between them rises from both ends to
their first common stabilization.  We call this route the \emph{canonical
component path}.  Its lift to triangulations is worth stating separately.

\begin{corollary*}[Canonical component paths]
Any two simplicial triangulations \(P\) and \(Q\) of the same closed
connected PL \(d\)-manifold admit a Pachner sequence of the form
\[
 P
 \xrightarrow{\;1\text{--}(d+1)\;*}
 P^+
 \xrightarrow{\;\text{vertex-preserving moves}\;*}
 Q^+
 \xrightarrow{\;(d+1)\text{--}1\;*}
 Q.
\]
Here and below, the asterisk denotes a finite sequence of moves of the
indicated type.  This sequence lifts the canonical component path
between the fixed-vertex components containing \(P\) and \(Q\).  The
number of vertices first increases, then remains constant, and finally
decreases.
\end{corollary*}

The component path is canonical, but a Pachner sequence lifting it is not: the
result does not select the facets to subdivide or a preferred horizontal flip
path.  This result is also distinct from the unimodal path theorem of Burton
and He \cite{BH26}.  Their theorem concerns generalized one-vertex
triangulations of \(3\)-manifolds and grades complexity by the number of
tetrahedra.  Here the triangulations are simplicial, the horizontal graphs fix
the number of vertices, and the component graph is graded by that number.
The hypotheses and grading are therefore different.

The component tree need not have a unique vertex at its lowest occupied
level.  Every component determines a unique stabilization ray, and any two
such rays eventually share a common tail.  A chosen triangulation therefore
determines a monotone closure under facet subdivisions and vertex-preserving
moves, but that choice need not be distinguished.

Let \(n_0\) be the minimum number of vertices in a triangulation of \(M\).
When \(G(n_0)\) is connected, its unique component determines a distinguished
stabilization ray, which we call the \emph{trunk}.  In particular, a manifold
has a trunk whenever its triangulation with the minimum number of vertices is
unique up to combinatorial isomorphism, as for a sphere or the \(2\)-torus.
The trunk is a ray of components; it does not select a preferred
triangulation, facet, or flip sequence at any level.  For \(S^d\), it begins
with the boundary of the \((d+1)\)-simplex.  The computations below describe
the first levels of the \(S^3\) trunk and several branches that merge into it.

The paper is organized as follows.  Section~\ref{sec:definitions} sets up
bistellar notation in arbitrary dimension and defines support, fixed-vertex
components, lift fibers, and unflippability.  Section~\ref{sec:stabilization}
proves the Component Stabilization Theorem.  Section~\ref{sec:tree} derives
the component tree, its predecessor structure, coalescence, and canonical
component paths as consequences.  Section~\ref{sec:rays} develops
stabilization rays, monotone closures, and trunks.  Section~\ref{sec:examples}
gives examples outside \(S^3\), including one-step mergers in surface flip
graphs, the K\"uhnel--Lassmann \(3\)-torus, and unflippable spheres in every
dimension at least four.  Sections~\ref{sec:s3} and \ref{sec:experiments}
specialize to triangulated \(3\)-spheres and present the computational results.

\section{Bistellar moves and fixed-vertex Pachner graphs}
\label{sec:definitions}

Throughout, \(M\) denotes a closed connected PL \(d\)-manifold with
\(d\ge2\), and a triangulation means a finite simplicial complex whose
geometric realization is PL-homeomorphic to \(M\).  We identify isomorphic
triangulations.

\paragraph{Basic simplicial terminology.}
A \emph{facet} of a \(d\)-dimensional triangulation is a \(d\)-simplex, and a
\emph{ridge} is a \((d-1)\)-simplex.  For a simplex \(\alpha\in T\), its
\emph{link} is
\[
 \operatorname{lk}_T(\alpha)
 =\{\gamma\in T:\gamma\cap\alpha=\varnothing
   \text{ and }\gamma\cup\alpha\in T\}.
\]
We write \(\partial\alpha\) for the boundary complex of \(\alpha\) and
\(A*B\) for the simplicial join of two complexes.  The closed star of
\(\alpha\) is \(\alpha*\operatorname{lk}_T(\alpha)\).

\paragraph{Bistellar moves.}
Let \(\alpha\) and \(\beta\) be disjoint simplices with
\[
 |\alpha|+|\beta|=d+2,
\]
where \(|\cdot|\) denotes the number of vertices.  Suppose that
\(\alpha\in T\), \(\beta\notin T\), and
\(\operatorname{lk}_T(\alpha)=\partial\beta\).  Thus, the closed star of
\(\alpha\) is the ball \(\alpha*\partial\beta\), and the bistellar move
\[
 \alpha*\partial\beta
 \quad\longleftrightarrow\quad
 \partial\alpha*\beta
\]
replaces this ball by the second one.  If \(|\beta|=p\) and
\(|\alpha|=q\), so that \(p+q=d+2\), we call the forward operation a
\(p\)--\(q\) move.  Its source contains \(p\) facets and its target contains
\(q\) facets.

The \(1\)--\((d+1)\) move is stellar subdivision of a facet and adds one
vertex.  Its inverse is the \((d+1)\)--\(1\) move and deletes one vertex.
All other bistellar moves preserve the vertex set; we call them
\emph{proper} or \emph{vertex-preserving} moves.  Equivalently, these are the
\(p\)--\(q\) moves with \(2\le p,q\le d\).

\paragraph{Support.}
The \emph{support} of a bistellar move in its source triangulation is the
bistellar ball \(\alpha*\partial\beta\), or, equivalently when only facets
are being counted, the collection of \(p\) facets removed by the move.

\paragraph{Pachner graphs and horizontal slices.}
Let \(\mathcal P\) be the full Pachner graph whose vertices are all
simplicial triangulations of the fixed manifold \(M\), with edges given by
all bistellar moves.  For each \(n\), let \(G(n)\) be the induced graph on
the triangulations with \(n\) vertices, with edges given by all proper
bistellar moves.  We retain a loop when such a move produces an isomorphic
triangulation.  We call \(G(n)\) the \emph{restricted Pachner graph} or
\emph{fixed-vertex flip graph}.  Let
\[
 C(n)=\pi_0\bigl(G(n)\bigr)
\]
be its set of connected components.  As in the introduction, the fixed
manifold is suppressed from the notation; when several manifolds are under
discussion, we write \(G_M(n)\) and \(C_M(n)\).

We write \(T_1\sim_n T_2\) when \(T_1\) and \(T_2\) lie in the same
component of \(G(n)\).

\paragraph{Lift fibers.}
If \(T\) is a vertex of \(G(n)\) and \(\sigma\) is a facet of \(T\), let
\(T^\sigma\) denote the triangulation obtained by stellarly subdividing
\(\sigma\) with a new vertex.  The \emph{lift fiber} of \(T\) is
\[
 L(T)=\{T^\sigma:\sigma\text{ is a facet of }T\}
 \subseteq V\bigl(G(n+1)\bigr).
\]

\paragraph{Unflippable triangulations.}
An \(n\)-vertex triangulation \(T\) of \(M\) is \emph{unflippable} if it
admits no proper bistellar move.  Equivalently, \(T\) is an isolated vertex
of \(G(n)\).

\paragraph{The component graph.}
The \emph{component graph} \(\mathcal C\) is the directed graded graph whose
vertices are the elements of \(C(n)\), over all \(n\).  There is a directed
edge from \(U\in C(n)\) to \(V\in C(n+1)\) if some facet subdivision of some
triangulation in \(U\) belongs to \(V\).

\section{Component stabilization in arbitrary dimension}
\label{sec:stabilization}

\begin{theorem}[Component Stabilization Theorem]
\label{thm:component-stabilization}
Every vertex of \(\mathcal C\) has out-degree one.  Equivalently, if
\(U\in C(n)\) and \(P,Q\in U\), then every element of \(L(P)\) and every
element of \(L(Q)\) lie in the same connected component of \(G(n+1)\).
\end{theorem}

We prove Theorem~\ref{thm:component-stabilization} through three lemmas.  The
first is a local transport move.

\begin{lemma}[Sliding a stellar vertex across a ridge]
\label{lem:slide}
Let \(T\) be a triangulation of a closed \(d\)-manifold.  Suppose the facets
\(a*F\) and \(b*F\) share the ridge \(F\).  Insert a new vertex \(v\) into
\(a*F\).  The resulting triangulation is connected by two proper bistellar
moves to the triangulation obtained by inserting \(v\) into \(b*F\).
More precisely, the slide is a \(2\)--\(d\) move followed by a
\(d\)--\(2\) move.
\end{lemma}

\begin{proof}
After subdividing \(a*F\), the two facets \(v*F\) and \(b*F\) form the
bistellar ball
\[
 F*\partial(vb).
\]
The edge \(vb\) is absent, so the \(2\)--\(d\) move
\[
 F*\partial(vb)\longrightarrow \partial F*(vb)
\]
is legal.  The ridge \(F\) is absent after this move.  The remaining facets
of the stellar subdivision of \(a*F\) form
\[
 \partial F*(va),
\]
so the \(d\)--\(2\) move
\[
 \partial F*(va)\longrightarrow F*\partial(va)
\]
is legal.  The second move restores \(a*F\), while the first move together
with the surviving facet \(v*F\) gives precisely the stellar subdivision of
\(b*F\).  Both moves are proper.  When \(d=2\), they are two diagonal flips.
\end{proof}

\begin{lemma}[Lift fibers are connected]
\label{lem:fibers-connected}
For every triangulation \(T\) of a closed connected \(d\)-manifold, all
elements of \(L(T)\) lie in one connected component of \(G(n+1)\).
\end{lemma}

\begin{proof}
The facet-adjacency graph of a connected triangulated closed manifold is
connected.  To see this, choose points in the interiors of two facets and
join them by a PL path in \(|T|\).  Since the \((d-2)\)-skeleton has
codimension at least two, PL general position allows the path to be chosen
disjoint from it.  The path therefore passes from one facet to another only
through the relative interiors of ridges, and the facets it traverses form a
path in the facet-adjacency graph.  Successive applications of
Lemma~\ref{lem:slide} transport the stellar vertex along such a path using
only proper bistellar moves.
\end{proof}

The third lemma lifts an edge of \(G(n)\) to an edge of \(G(n+1)\).

\begin{lemma}[Proper edges lift between fibers]
\label{lem:edges-lift}
If \(P\) and \(Q\) are adjacent in \(G(n)\), then there are
\(P'\in L(P)\) and \(Q'\in L(Q)\) that are adjacent in \(G(n+1)\).
\end{lemma}

\begin{proof}
Suppose \(P\to Q\) is a proper \(p\)--\(q\) move, written
\[
 \alpha*\partial\beta\longrightarrow\partial\alpha*\beta.
\]
Its support contains \(p\le d\) facets.  Each facet of a closed simplicial
\(d\)-manifold has \(d+1\) ridges, and every ridge belongs to exactly two
facets.  Two distinct facets share at most one ridge, because their
intersection is either empty or a single common face of dimension at most
\(d-1\).  Hence the \(d+1\) ridges of a facet lead to \(d+1\) distinct
neighbors, so the facet-adjacency graph is simple and \((d+1)\)-regular.  It
therefore has at least \(d+2\) vertices.  Thus, \(P\) has a facet \(\rho\)
outside the support.

The move leaves \(\rho\) and every proper face of \(\rho\) unchanged.
Stellar subdivision of \(\rho\) leaves its boundary unchanged; among faces
on the original vertices it removes only \(\rho\), and every new face
contains the new vertex \(v\).  The complementary simplex \(\beta\) remains
absent as well: it uses only the original vertices, so subdivision cannot
create it, and \(\beta\notin P\) implies \(\beta\notin P^\rho\).
Consequently, the original bistellar ball and its attaching boundary are
unchanged and the complementary simplex is still absent, so the same proper
move is legal after subdivision.  Subdivide \(\rho\) in both \(P\) and
\(Q\).  The original move then gives an edge
\[
 P^\rho\longleftrightarrow Q^\rho
\]
in \(G(n+1)\).
\end{proof}

\begin{proof}[Proof of Theorem~\ref{thm:component-stabilization}]
Let \(U\in C(n)\).  Choose triangulations \(P,Q\in U\) and a path
\[
 P=S_0,S_1,\ldots,S_m=Q
\]
in \(G(n)\).  Lemma~\ref{lem:edges-lift} supplies, for each \(i\), an edge
between one element of \(L(S_i)\) and one element of \(L(S_{i+1})\).
Lemma~\ref{lem:fibers-connected} connects any two choices inside each fiber.
Concatenating the fiber paths and lifted edges connects any prescribed
element of \(L(P)\) to any prescribed element of \(L(Q)\).  Thus, all facet
subdivisions of all triangulations in \(U\) lie in one component of
\(G(n+1)\), and \(U\) has exactly one outgoing neighbor in \(\mathcal C\).
\end{proof}

For \(U\in C(n)\), denote its unique outgoing neighbor by \(s_n(U)\).  Thus,
Theorem~\ref{thm:component-stabilization} defines a map
\[
 s_n:C(n)\longrightarrow C(n+1).
\]

\begin{remark}[A quantitative lift bound]
The proof is constructive.  If a path
\(S_0,S_1,\ldots,S_m\) is given downstairs and \(D_i\) is the diameter of
the facet-adjacency graph of \(S_i\), then each base move can be lifted after
at most \(D_i\) slides, each slide costing two proper moves.  Including a
final slide to a prescribed terminal facet subdivision gives the upper bound
\[
 \sum_{i=0}^{m-1}(2D_i+1)+2D_m
\]
on the length of the lifted path.
\end{remark}

\section{Consequences: the component tree and canonical component paths}
\label{sec:tree}

All structural statements in this section are consequences of
Theorem~\ref{thm:component-stabilization} together with Pachner's theorem.
The underlying simple graph of \(\mathcal C\) is obtained from the full
Pachner graph \(\mathcal P\) by identifying all triangulations in each
component of \(G(n)\), and then suppressing loops and parallel edges.

\begin{corollary}[The component graph is a tree]
\label{cor:component-tree}
For every closed connected PL manifold \(M\) of dimension at least two, the
underlying undirected graph of \(\mathcal C\) is a tree.
\end{corollary}

\begin{proof}
By Pachner's theorem, the full Pachner graph \(\mathcal P\) is connected, so
its quotient obtained by identifying horizontal components is connected.
Suppose that the underlying undirected graph of \(\mathcal C\) contained a
cycle.  A vertex of minimum level on that cycle would have two distinct
neighbors at the next level, contradicting
Theorem~\ref{thm:component-stabilization}.  Thus, the underlying graph is
acyclic and, being connected, is a tree.
\end{proof}

\paragraph{Newborn components.}
A component of a fixed-vertex graph is \emph{newborn} if it has no neighbor at
the preceding level of \(\mathcal C\).

\begin{proposition}[Predecessors and leaves]
\label{prop:predecessors}
Let \(U\in C(n)\).  Then \(U\) has a neighbor at level \(n-1\) if and only
if some triangulation in \(U\) admits a \((d+1)\)--\(1\) bistellar move.
Consequently, \(U\) is newborn if and only if none of its triangulations
admits a vertex-deleting move, and this holds if and only if \(U\) is a leaf
of the component tree.
\end{proposition}

\begin{proof}
If \(U=s_{n-1}(V)\), choose a triangulation \(P\in V\) and a facet \(\rho\)
of \(P\).  The facet subdivision \(P^\rho\) belongs to \(U\), and its
stellar vertex supports the inverse \((d+1)\)--\(1\) move.

Conversely, suppose that \(T\in U\) admits a \((d+1)\)--\(1\) move.  The
move produces an \((n-1)\)-vertex triangulation \(P\), and \(T\) is a facet
subdivision of \(P\).  If \(V\) is the component containing \(P\), then
\(U=s_{n-1}(V)\).

Every node of \(\mathcal C\) has exactly one neighbor at the next level, and
its neighbors at the preceding level are precisely its predecessors.  It is
therefore a leaf exactly when it has no predecessor.
\end{proof}

\begin{corollary}[Unflippable leaves]
\label{cor:unflippable-leaves}
Every unflippable triangulation of \(M\) represents a newborn leaf of the
component tree.
\end{corollary}

\begin{proof}
Let \(T\) be unflippable.  Suppose that a \((d+1)\)--\(1\) move deletes a
vertex \(v\), so that \(\operatorname{lk}_T(v)=\partial\sigma\), where the
\(d\)-simplex \(\sigma\) is absent from \(T\).  Choose a ridge \(\tau\) of
\(\sigma\).  Since \(T\) is closed, the ridge \(\tau\) lies in exactly two
facets, \(v*\tau\) and \(w*\tau\).  The absence of \(\sigma\) implies
\(w\notin\sigma\), while \(\operatorname{lk}_T(v)=\partial\sigma\) implies
that \(vw\notin T\).  Thus,
\[
 \operatorname{lk}_T(\tau)=\partial(vw),
\]
so \(\tau*\partial(vw)\) supports a proper \(2\)--\(d\) bistellar move,
contrary to unflippability.  Thus, \(T\) admits no vertex-deleting move.
Since its fixed-vertex component consists only of \(T\),
Proposition~\ref{prop:predecessors} applies.
\end{proof}

Thus, stabilization may identify distinct components, but it cannot split one
component into several.  More strongly, any two stabilization rays merge.

\begin{corollary}[Coalescence]
\label{cor:coalescence}
Let \(U\in C(n)\) and \(V\in C(m)\).  There are integers \(r,s\ge0\)
such that the \(r\)-fold stabilization of \(U\) equals the \(s\)-fold
stabilization of \(V\).  In particular, any two components have a
common stabilization.
\end{corollary}

\begin{proof}
Take the unique finite path from \(U\) to \(V\) in the tree underlying
\(\mathcal C\).  Along this path there can be no strict local minimum of the
level function, since such a vertex would have two distinct upward neighbors.
Therefore, the path rises monotonically from both ends to one highest vertex.
That vertex is an iterated stabilization of both \(U\) and \(V\).
\end{proof}

\begin{corollary}[Canonical component paths]
\label{cor:canonical-path}
Let \(U\in C(n)\) and \(V\in C(m)\).  Their stabilization rays have a first
common component \(W\).  The unique path from \(U\) to \(V\) in the
component tree rises from \(U\) and \(V\) to \(W\); this is the
\emph{canonical component path} between \(U\) and \(V\).

If \(P\in U\) and \(Q\in V\), then this component path lifts to a Pachner
sequence
\[
 P
 \xrightarrow{\;1\text{--}(d+1)\;*}
 P^+
 \xrightarrow{\;\text{proper moves}\;*}
 Q^+
 \xrightarrow{\;(d+1)\text{--}1\;*}
 Q,
\]
where \(P^+\) and \(Q^+\) lie in \(W\).  Thus, the number of vertices first
increases, then remains constant, and finally decreases.
\end{corollary}

\begin{proof}
Corollary~\ref{cor:coalescence} gives a common component of the two
stabilization rays.  Since each ray has a unique successor at every level,
their common components form a common tail and therefore have a first member
\(W\).  The unique path in the tree rises from both ends to \(W\).

Choose any sequence of facet subdivisions of \(P\) and of \(Q\) that reaches
the level of \(W\).  Iterated applications of
Theorem~\ref{thm:component-stabilization} place the resulting triangulations
\(P^+\) and \(Q^+\) in \(W\).  By the definition of a component, they are
connected by proper bistellar moves.  Reversing the subdivisions of \(Q\)
gives the final descending segment.
\end{proof}

\begin{remark}
The component path is canonical; a Pachner sequence lifting it is not.  The
result does not select the facets to subdivide or a preferred proper-move path
inside the common component.  It also gives no useful a priori bound on the
height of the first common stabilization; obtaining effective bounds is one
of the main questions raised below.
\end{remark}

\section{Stabilization rays, monotone closures, and trunks}
\label{sec:rays}

For \(U\in C(n_0)\), the sequence
\[
 U,\quad s_{n_0}(U),\quad
 s_{n_0+1}s_{n_0}(U),\quad\ldots
\]
is the \emph{stabilization ray} from \(U\).  Corollary~\ref{cor:coalescence}
says that any two stabilization rays eventually share a common tail.  The
component tree is therefore one-ended, but it need not have a preferred
lowest vertex.

\paragraph{Monotone closure.}
Let \(R\) be a triangulation of \(M\) with \(n_0\) vertices.  Its
\emph{monotone closure} \(\mathcal M(R)\) is the set of triangulations
reachable from \(R\) by proper bistellar moves and
\(1\)--\((d+1)\) facet subdivisions, without using a
\((d+1)\)--\(1\) move.  For \(n\ge n_0\), let
\[
 \mathcal M_R(n)
 =G(n)\bigl[\mathcal M(R)\cap V(G(n))\bigr]
\]
be the induced subgraph at level \(n\).

\begin{corollary}[Monotone closures]
\label{cor:monotone-closure}
For every \(n\ge n_0\), the graph \(\mathcal M_R(n)\) is exactly one
connected component of \(G(n)\).  Its component is the
\((n-n_0)\)-fold stabilization of the component containing \(R\).
\end{corollary}

\begin{proof}
At level \(n_0\), the monotone closure contains exactly the component of
\(G(n_0)\) containing \(R\).  Suppose its level-\(n\) slice is one component
\(U\).  Every monotone path to level \(n+1\) has a last facet subdivision,
followed only by proper moves.  By the Component Stabilization Theorem, all
facet subdivisions of all triangulations in \(U\) lie in the single
component \(s_n(U)\), and proper moves fill out that entire component.  Thus,
the next slice is exactly \(s_n(U)\).  Induction proves the claim.
\end{proof}

In particular, the first move in a Pachner path that starts in
\(\mathcal M(R)\) and leaves \(\mathcal M(R)\) must be a
\((d+1)\)--\(1\) move: proper moves and facet subdivisions remain in the
monotone closure.

\begin{corollary}[Eventual merger with a chosen ray]
\label{cor:absorption}
For every component \(U\in C(n)\) and every triangulation \(R\), some
iterated stabilization of \(U\) is a level of \(\mathcal M(R)\).
\end{corollary}

\begin{proof}
Apply Corollary~\ref{cor:coalescence} to \(U\) and the component containing
\(R\).  The upward orbit of the latter component is precisely the sequence of
components in the monotone closure.
\end{proof}

\paragraph{Trunks.}
Let \(n_0\) be the minimum number of vertices in a triangulation of \(M\).
Suppose that \(G(n_0)\) is connected, and let \(U_0\) be its unique
component.  The stabilization ray beginning at \(U_0\) is the \emph{trunk}
of \(M\).  We write \(\mathcal T(n)\) for its level-\(n\) component.  By
Corollary~\ref{cor:monotone-closure}, the triangulations in these components
form \(\mathcal M(R)\) for any \(R\in U_0\).

In particular, \(M\) has a trunk whenever its triangulation with the minimum
number of vertices is unique up to combinatorial isomorphism.  If the graph
at the minimum level has more than one component, there is no distinguished
trunk, although every component determines a stabilization ray and all of
those rays eventually have a common tail.  A trunk is a ray of components; it
does not select a preferred triangulation, facet, or flip sequence at any
level.

\paragraph{The sphere trunk.}
For \(S^d\), the boundary \(\partial\Delta^{d+1}\) of the
\((d+1)\)-simplex is the unique triangulation with \(d+2\) vertices and hence
begins the trunk.  For \(d=3\), let
\[
 \mathcal T(n)=\mathcal M_{\partial\Delta^4}(n).
\]

\begin{corollary}[Sphere trunk components]
\label{cor:trunk-component}
For every \(n\ge5\), \(\mathcal T(n)\) is exactly one connected component
of \(G(n)\).
\end{corollary}

\section{Examples beyond triangulated \texorpdfstring{\(S^3\)}{S3}}
\label{sec:examples}

The theorem and its consequences apply to every closed connected PL manifold
of dimension at least two.  The following examples illustrate both the
presence of several components at the minimum level and the speed with which
distinct stabilization rays can merge.

\subsection{Surface flip graphs and one-step coalescence}
\label{sec:surfaces}

Because several manifolds occur in this subsection, we restore the suppressed
subscript and write \(G_F(n)\) and \(C_F(n)\) for a surface \(F\).
In dimension two, the proper bistellar move is the \(2\)--\(2\) diagonal
flip, and facet subdivision is the \(1\)--\(3\) move.  If a closed
triangulated surface has \(n\) vertices, Euler characteristic \(\chi\), and
\(e\) edges, then
\[
 e=3(n-\chi).
\]
Equality \(e=\binom n2\) forces the triangulation to be neighborly, and a
neighborly surface triangulation admits no diagonal flip, because the
opposite diagonal is already present.  Following Negami
\cite{N99frozen}, a surface triangulation admitting no diagonal flip is
called \emph{frozen}; frozen is the two-dimensional case of unflippability
as defined in Section~\ref{sec:definitions}.

Fixed-vertex flip connectivity of surfaces has a long history.  Wagner proved
connectivity for the sphere, and Dewdney treated the torus \cite{W36,D73}.
More generally, Negami proved that for every closed surface \(F\), the graph
\(G_F(n)\) is connected for all sufficiently large \(n\)
\cite{N94,N99survey}; King later gave an explicit bound linear in
\(\lvert\chi(F)\rvert\) \cite{King03}.  A related result of Komuro,
Nakamoto, and Negami preserves minimum degree at least four \cite{KNN99}.
Thus, the component tree of a surface eventually has a single ray.  Our
computations instead concern its first mergers at minimum levels.

The unique \(7\)-vertex triangulation of the torus is the classical
M\"obius triangulation, whose \(1\)-skeleton is \(K_7\); its polyhedral
realization is due to Cs\'asz\'ar \cite{Csaszar49,K23}.  It is frozen and
begins the trunk of \(T^2\).  The following four minimum levels have several
components.  Here \(N_h\) denotes the nonorientable surface of genus \(h\),
and \(\Sigma_g\) the orientable surface of genus \(g\).  The minimum numbers of vertices follow from Ringel and
Jungerman--Ringel \cite{R55,JR80}; the enumerations and component decompositions use the surface censuses
\cite{Sul06,SulPM,SL09,LutzSurf}, including the classification of the \(59\)
neighborly \(\Sigma_6\) triangulations \cite{ABS96}.

\begin{center}
\begin{tabular}{c|c|r|c}
surface & minimum level & triangulations & component sizes \\ \hline
\(N_4\)       & 9  & 37 & \(32,3,2\) \\
\(N_5\)       & 9  & 2  & \(1,1\) \\
\(\Sigma_3\)  & 10 & 20 & \(13,5,1,1\) \\
\(\Sigma_6\)  & 12 & 59 & \(1^{59}\)
\end{tabular}
\end{center}

The two triangulations of \(N_5\) and all \(59\) triangulations of
\(\Sigma_6\) are neighborly and therefore frozen.  The nontrivial component
sizes for \(N_4\) and \(\Sigma_3\) show that minimum-level branches need not
be singleton components.

\begin{proposition}[One-step coalescence at minimum levels]
\label{prop:surface-coalescence}
The stabilization maps
\[
 s_9^{N_4}:C_{N_4}(9)\to C_{N_4}(10),\qquad
 s_9^{N_5}:C_{N_5}(9)\to C_{N_5}(10),
\]
\[
 s_{10}^{\Sigma_3}:C_{\Sigma_3}(10)\to C_{\Sigma_3}(11),\qquad
 s_{12}^{\Sigma_6}:C_{\Sigma_6}(12)\to C_{\Sigma_6}(13)
\]
are constant: in each case, all components at the minimum level merge into
one component after a single facet subdivision.  For \(N_4\)
and \(\Sigma_3\), the maps are moreover surjective, since
\(G_{N_4}(10)\) and \(G_{\Sigma_3}(11)\) are
connected.  In particular, a single stabilization fiber can contain at least
\(59\) fixed-vertex components.
\end{proposition}

\begin{proof}
We used the catalogued facet lists from the Manifold Page \cite{LutzSurf}.
For \(N_4\), the complete catalog of \(13{,}657\) ten-vertex triangulations
forms one component; the same is true of the \(65{,}878\) eleven-vertex
triangulations of \(\Sigma_3\).  This proves both constancy and surjectivity
in those two cases.

For \(N_5\) and \(\Sigma_6\), explicit diagonal-flip paths join a chosen
facet subdivision of one representative from every source component.  The
Component Stabilization Theorem makes the representative and facet choices
immaterial, so both maps are constant.  These certificates do not enumerate
the entire next level and therefore do not assert that
\(G_{N_5}(10)\) or \(G_{\Sigma_6}(13)\) is connected.
\end{proof}

\subsection{Unflippable triangulations beyond \texorpdfstring{\(S^3\)}{S3}}

The classical \(15\)-vertex triangulation of the \(3\)-torus constructed by
K\"uhnel and Lassmann is neighborly, and its edges have valence \(4\) or
\(6\) \cite{KL84,BK12}.  Neighborliness blocks every \(2\)--\(3\) move,
and the absence of a valence-three edge blocks every \(3\)--\(2\) move.
It is therefore unflippable.  It is often described as the smallest known
triangulation of \(T^3\), although vertex-minimality and uniqueness remain
open \cite{L09,BK12}.

Spreer's census of cyclic combinatorial \(3\)-manifolds contains \(6070\)
combinatorial types with at most \(22\) vertices, representing \(67\)
topological types, and records locally minimal examples at many levels
\cite{S14}.  Here ``locally minimal'' means that the triangulation cannot be
reduced by bistellar moves without first adding a vertex; this is weaker than
unflippability, but the census supplies a broad source of non-spherical
examples for the fixed-vertex questions considered here.  The known
unflippable spheres \(U(20),U_1(21),U_2(21)\) also arise from this census.

The Dougherty--Faber--Murphy sphere \(U(16)\) is an unflippable
\(3\)-sphere \cite{DFM04}.  Bagchi and Datta use it to construct
unflippable spheres in all larger dimensions \cite{BD13}.  They first build a
\(4\)-ball \(B^4_{16}\) with boundary \(U(16)\).  If
\(B^e_{e+1}\) is the standard \(e\)-simplex, then for every \(e\ge0\),
\[
 S^{e+4}_{e+17}
 =\partial\bigl(B^4_{16}*B^e_{e+1}\bigr)
\]
is an unflippable combinatorial \((e+4)\)-sphere.  In dimension \(d=e+4\),
this gives an unflippable \(d\)-sphere on \(d+13\) vertices for every
\(d\ge4\).

Each such sphere is an isolated fixed-vertex component.  Nevertheless, Corollary~\ref{cor:coalescence} guarantees that its stabilization
ray eventually
merges with the sphere trunk.  Thus, unflippable triangulations
appear in arbitrarily high dimensions, while stabilization imposes a uniform
global constraint on how their components fit into the full Pachner graph.

\section{The trunk of triangulated \texorpdfstring{\(S^3\)}{S3}}
\label{sec:s3}

In this section and Section~\ref{sec:experiments}, fix \(M=S^3\).  Thus,
\(G(n)\) and \(C(n)\) denote the fixed-vertex graph and its component set for
triangulated \(3\)-spheres.  The monotone closure of \(\partial\Delta^4\)
consists of all triangulations reachable using
\(1\)--\(4\), \(2\)--\(3\), and \(3\)--\(2\) moves, but no
\(4\)--\(1\) moves.  By Corollary~\ref{cor:trunk-component}, its
level-\(n\) slice \(\mathcal T(n)\) is a single induced component of
\(G(n)\).

\paragraph{Polytopal spheres lie in the trunk.}
Every simplicial \(4\)-polytope admits a Schlegel diagram, which gives a
regular triangulation of a point configuration in \(\mathbb R^3\) whose
convex hull is a tetrahedron.  By the incremental topological flipping
theorem of Edelsbrunner and Shah \cite{ES96}, such a regular triangulation
can be constructed from the outer tetrahedron by inserting points and
applying \(2\)--\(3\) and \(3\)--\(2\) flips.  In the full Pachner graph,
this gives a path from \(\partial\Delta^4\) to the boundary of the given
polytope using only \(1\)--\(4\), \(2\)--\(3\), and \(3\)--\(2\) moves.
Thus, every polytopal triangulated \(3\)-sphere with \(n\) vertices lies in
\(\mathcal T(n)\).

\paragraph{The trunk also contains non-polytopal spheres.}
Non-polytopal triangulated \(3\)-spheres first occur with eight vertices.  The
Br\"uckner and Barnette spheres are the two classical examples
\cite{GS67,Bar73}.  A direct check shows that both lie in
\(\mathcal T(8)\): each is one \(2\)--\(3\) flip from a polytopal
triangulation or from the other, so \(G(8)=\mathcal T(8)\).

At \(n=9\), Altshuler, Bokowski, and Steinberg \cite{ABS80} showed that
\(G(9)\) is connected, so \(G(9)=\mathcal T(9)\).  Their classification
contains \(1296\) simplicial \(3\)-spheres.  Using Firsching's classification
and data \cite{F20,F-Data}, we identified the \(1142\) polytopal and \(154\)
non-polytopal types and verified that the subgraph induced by the
non-polytopal types is connected.

\paragraph{Known unflippable 3-spheres.}
The four currently known examples relevant to our computations are
\[
 U(16),\qquad U(20),\qquad U_1(21),\qquad U_2(21).
\]
The \(16\)-vertex example is from \cite{DFM04}; the \(20\)- and
\(21\)-vertex examples arise among the cyclic triangulations studied in
\cite{S14}.  Each is neighborly and has minimum edge valence at least four.
Thus, every \(2\)--\(3\) move is blocked by an existing complementary edge
and every \(3\)--\(2\) move is blocked by edge valence.  Each is therefore
unflippable and, by Corollary~\ref{cor:unflippable-leaves}, represents a
newborn singleton leaf of the component tree, distinct from
\(\mathcal T(n)\) for the corresponding value of \(n\).

The computational results below extend the connected range through
\(n=11\), treat a high-valence part of the \(n=12\) seed space, and show
that these four isolated branches merge into the trunk
after one stabilization.

\section{Computational results for \texorpdfstring{\(S^3\)}{S3}}
\label{sec:experiments}

The structural results above reduce membership in \(\mathcal T(n)\) to one
fixed-vertex connectivity problem.  Throughout this section, we use the
computational framework of Sulanke and Lutz \cite{SL09} for seed enumeration
and, for the \(10\)-vertex polytopality labels, Firsching's classification
and accompanying online data \cite{F17,F-Data}.

\subsection{\texorpdfstring{\(G(10)\)}{G(10)} is connected}

Lutz enumerated the \(247{,}882\) ten-vertex triangulated
\(3\)-spheres \cite{Lutz08}.  Firsching classified these as \(162{,}004\)
polytopal and \(85{,}878\) non-polytopal types \cite{F17,F-Data}.  We use the
Lutz census identifiers, as recorded in Firsching's data, and use the
polytopal labels only as a convenient source of vertices already known to lie
in \(\mathcal T(10)\).

\paragraph{Seed triangulations.}
A triangulation is a \emph{seed triangulation} if it admits no legal
\(3\)--\(2\) bistellar flip.

Every triangulation is either a seed or can be connected to a seed by a
sequence of \(3\)--\(2\) flips.  Indeed, each legal \(3\)--\(2\) move
reduces the number of tetrahedra, so repeated \(3\)--\(2\) moves terminate.

Thus, to prove that all of \(G(10)\) lies in \(\mathcal T(10)\), it
suffices to check the ten-vertex seeds.  Taking the complement of Firsching's
polytopal list and testing each non-polytopal type for legal \(3\)--\(2\)
moves, we found exactly \(15\) non-polytopal seeds.  For each, we enumerated
its legal \(2\)--\(3\) moves and found a move to a polytopal triangulation.
Since polytopal triangulations lie in \(\mathcal T(10)\), each of the \(15\)
non-polytopal seeds also lies in \(\mathcal T(10)\).

The remaining seeds are polytopal and hence also lie in \(\mathcal T(10)\).
Therefore every \(10\)-vertex seed lies in \(\mathcal T(10)\), and every
\(10\)-vertex triangulated \(3\)-sphere is connected by \(3\)--\(2\) flips
to \(\mathcal T(10)\).  Consequently,
\[
 G(10)=\mathcal T(10),
\]
and \(G(10)\) is connected.

\subsection{\texorpdfstring{\(G(11)\)}{G(11)} is connected}

Using Sulanke and Lutz's software \cite{SL09}, we enumerated the \(681\) seed
triangulations with \(11\) vertices.  After preliminary topological filters,
\(225\) candidate \(3\)-manifold seeds remained.  Simulated annealing and
isomorphism testing identified \(222\) as triangulated \(3\)-spheres; the
remaining three were triangulations of \(\mathbb{RP}^3\).

Since \(G(10)\) is connected, every \(10\)-vertex triangulated
\(3\)-sphere lies in \(\mathcal T(10)\).  Hence every \(1\)--\(4\)
subdivision of a \(10\)-vertex sphere lies in the single stabilized component
\(\mathcal T(11)\).  To prove that \(G(11)\) is connected, it is
enough to connect every \(11\)-vertex seed sphere to one such subdivision.

Starting from \(1\)--\(4\) subdivisions of \(10\)-vertex spheres, we
performed random walks in \(G(11)\).  These walks produced explicit
restricted-flip paths from vertices of \(\mathcal T(11)\) to all \(222\) spherical
\(11\)-vertex seeds.  Therefore every \(11\)-vertex seed sphere lies in
\(\mathcal T(11)\).

Every non-seed admits a \(3\)--\(2\) move, and repeated \(3\)--\(2\) moves
terminate at a seed.  Reversing these moves stays within \(G(11)\).
Hence
\[
 G(11)=\mathcal T(11),
\]
and \(G(11)\) is connected.

\subsection{A high-valence seed computation in \texorpdfstring{\(G(12)\)}{G(12)}}

Full seed enumeration at \(n=12\) is substantially more expensive than at
\(n=11\).  The most plausible candidates for new components outside
\(\mathcal T(12)\) are high-valence seeds, since the known unflippable examples have no
valence-three edges.  We therefore focused first on the high-valence regime.

We exhaustively enumerated all \(12\)-vertex seed triangulations whose edges
all have valence at least \(4\).  For every seed in this class, we found an
explicit restricted-flip path to \(\mathcal T(12)\).  Equivalently, the union
of these paths gives a connected subgraph of \(G(12)\) containing
all high-valence seeds and a known vertex of \(\mathcal T(12)\).  Thus, no
high-valence \(12\)-vertex seed gives a component outside
\(\mathcal T(12)\).

Any possible disconnected component of \(G(12)\) must therefore
come from a seed containing a valence-three edge whose \(3\)--\(2\) move is
blocked by the complementary triangle.  A complete proof that
\(G(12)\) is connected requires enumerating and analyzing these
blocked valence-three seeds.  For larger \(n\), exhaustive seed enumeration
becomes rapidly more difficult because of the growth rate of triangulated
\(3\)-spheres \cite{NSW16}.

\subsection{The known unflippable spheres after one facet subdivision}

For each of the four unflippable spheres
\[
 U(16),\quad U(20),\quad U_1(21),\quad U_2(21),
\]
we performed one \(1\)--\(4\) move to obtain a triangulation \(U^+\) on one
additional vertex.  We then found an explicit simplification sequence from
\(U^+\) to \(\partial\Delta^4\) using only \(2\)--\(3\), \(3\)--\(2\),
and \(4\)--\(1\) moves.  Reversing the sequence gives a path from
\(\partial\Delta^4\) to \(U^+\) using only \(1\)--\(4\), \(2\)--\(3\),
and \(3\)--\(2\) moves.

\begin{proposition}[One-step merger of the known unflippable spheres]
\label{prop:one-step-absorption}
For each of
\[
 U(16),\quad U(20),\quad U_1(21),\quad U_2(21),
\]
every \(1\)--\(4\) subdivision lies in \(\mathcal T(n+1)\), where
the original sphere has \(n\) vertices.  Thus, each isolated component
merges into the sphere trunk after one stabilization.
\end{proposition}

\begin{proof}
For each example, the accompanying data contains an explicit sequence from
one subdivision \(U^+\) to \(\partial\Delta^4\) using only \(2\)--\(3\),
\(3\)--\(2\), and \(4\)--\(1\) moves.  Reversing the sequence shows that
\(U^+\in\mathcal T(n+1)\).  Lemma~\ref{lem:fibers-connected} then gives the
same conclusion for every facet subdivision of \(U\).
\end{proof}

The simplification sequences were found by simulated annealing
\cite{BL00,JLL22}, but the annealing procedure is not part of the proof.  Once
a sequence is found, it is a directly checkable certificate.  These examples
show that the coalescence guaranteed by Corollary~\ref{cor:coalescence}
can occur after only one stabilization.

\begin{corollary}[First branching level for \(S^3\)]
\label{cor:first-branching-bound}
Let
\[
 b(S^3)=\min\{n:G(n)\text{ is disconnected}\}.
\]
Then
\[
 12\le b(S^3)\le16.
\]
\end{corollary}

\smallskip
\begin{proof}
The fixed-vertex graphs are connected through level \(11\), so
\(b(S^3)\ge12\).  The unflippable sphere \(U(16)\) represents a newborn leaf
by Corollary~\ref{cor:unflippable-leaves}, whereas the level-\(16\) trunk
component has predecessors.  They are therefore distinct components of
\(G(16)\), and \(b(S^3)\le16\).
\end{proof}

\paragraph{Data availability.}
The computational certificates used in this paper are available in the
archival Zenodo record \url{https://doi.org/10.5281/zenodo.21980010}.  These
include the seed lists for the \(n=10\),
\(n=11\), and high-valence \(n=12\) computations; the restricted-flip paths
certifying membership in the relevant \(\mathcal T(n)\); the simplification
sequences for \(U(16)^{+}\), \(U(20)^{+}\), \(U_1(21)^{+}\), and \(U_2(21)^{+}\);
and the code and connecting paths for the surface stabilization experiments
of Proposition~\ref{prop:surface-coalescence}.  Each certificate can be
checked independently by verifying every listed move and the claimed terminal
triangulation.

\section{Conclusions and questions}
\label{sec:conclusions}

Return to an arbitrary fixed closed connected PL manifold \(M\).  The main
theorem gives a global organization of its Pachner graph.  Fixed-vertex
components form a one-ended graded tree under facet stabilization; any two rays eventually merge; and every pair of
triangulations has a canonical component path with a three-stage Pachner
lift.  For \(S^3\), the boundary of the \(4\)-simplex begins the trunk.  Its level components
contain every triangulated \(3\)-sphere through \(11\) vertices and merge
with each of the four known isolated sphere components after one facet
subdivision.

The component-tree formulation suggests several natural open questions.
Eventual absorption is automatic; the central issues are now
quantitative bounds, the geometry of branching, and the first levels at which
new branches appear.

\begin{enumerate}

\item[\textbf{(Q1)}] \emph{Effective stabilization and absorption bounds.}
Given triangulations \(P\) and \(Q\) of a fixed \(d\)-manifold, bound the
least level at which iterated facet subdivisions of \(P\) and \(Q\) enter
one restricted component.  Can this height be bounded in terms of the
numbers of vertices or facets, the length of a known Pachner path, or other
local data?  For the \(S^3\) trunk and \(U\in C(n)\), let \(a(U)\)
be the least \(k\ge0\) for which the \(k\)-fold stabilization of \(U\) is
\(\mathcal T(n+k)\).  The four known unflippable spheres have \(a(U)=1\).
How large can \(a(U)\) be?  General
simplification bounds such as those of Mijatovi\'c \cite{M03} are far too
large to describe the observed behavior.

\item[\textbf{(Q2)}] \emph{Branching and merging.}
Characterize the fibers of
\[
 s_n:C(n)\longrightarrow C(n+1).
\]
The surface computations show that a single stabilization fiber can contain
at least \(59\) components: all \(59\) isolated components at the minimum
level for \(\Sigma_6\) merge in one step.  More generally, is the stabilization map
from the minimum level of every closed surface constant?  Are the fiber
sizes unbounded, even for surfaces?  When complete target catalogs are
available, does the common image component exhaust the next fixed-vertex
graph?  More generally, can a newborn component of a surface occur above the
minimum level?

For \(S^3\), Corollary~\ref{cor:first-branching-bound} gives
\(12\le b(S^3)\le16\).  What is the exact first branching level?  At
\(n=12\), only seeds with blocked valence-three edges remain to be analyzed.
Does a component
outside \(\mathcal T(n)\) contain more than one triangulation?  The
\(N_4\) surface example shows that non-singleton branches occur naturally,
although every presently known component outside the sphere trunk is a
singleton.

\item[\textbf{(Q3)}] \emph{Geometry of branches outside the sphere trunk.}
All four currently known unflippable triangulated \(3\)-spheres are
neighborly.  Does there exist a non-neighborly unflippable triangulated
\(3\)-sphere?  More broadly, does every
triangulated \(3\)-sphere admitting a \emph{\(3\)-diagram} lie in
\(\mathcal T(n)\) when it has \(n\) vertices?  Here a \(3\)-diagram
means a choice of an exterior facet whose complement has a straight-line
realization in \(\mathbb R^3\) with tetrahedral convex hull.  Regular \(3\)-diagrams lie in the sphere trunk by
the theorem of Edelsbrunner and Shah.  A positive answer
would separate geometric
obstructions to the flip-connectivity problem of Edelsbrunner, Preparata,
and West \cite{EPW90} from the abstract component structure.

\item[\textbf{(Q4)}] \emph{Higher-dimensional component trees.}
Bagchi and Datta provide unflippable leaves in every dimension at
least four, while recent computations show that simplification in dimension
four can encounter substantial combinatorial depth \cite{BBS26}.  What are
the first branching levels for \(S^d\), how complicated can the horizontal
components and stabilization fibers become, and how rapidly do branches
merge into the sphere trunk?  Corollary~\ref{cor:component-tree} gives
stable connectivity but no practical strategy for finding the horizontal
path at a common stabilization level.

\end{enumerate}

\section*{Acknowledgements}

The authors thank OpenAI's ChatGPT for assistance during an
exploratory phase of the proof of
Theorem~\ref{thm:component-stabilization}.  We thank Anthropic's
Claude for coding assistance during exploratory phases of this work,
including the surface computations of
Proposition~\ref{prop:surface-coalescence}.  All proofs, computations,
and references were checked and edited by the authors, who take sole
responsibility for the final text.

\end{document}